\documentstyle[12pt]{article}

\def\be{\begin{equation}}
\def\ee{\end{equation}}
\def\bea{\begin{eqnarray}}
\def\eea{\end{eqnarray}}
\begin{document}

\title{\large{\bf GEOMETRY OF QUANTUM HOMOGENEOUS SUPERVECTOR BUNDLES 
AND REPRESENTATION THEORY OF QUANTUM GENERAL LINEAR SUPERGROUP}} 

\author{ R. B. ZHANG\\ 
 Department of Pure Mathematics\\
 University of Adelaide\\
 North Terrace, Adelaide\\
 S.A. 5001, Australia}
\date{ }
\maketitle

\vspace{1.5cm}

The quantum general linear supergroup $GL_q(m|n)$ 
is defined and its structure is studied systematically. 
Quantum homogeneous supervector bundles are introduced 
following Connes' theory, and applied to develop
the representation theory of $GL_q(m|n)$. 
Quantum Frobenius reciprocity is proven, and a  
Borel - Weil theorem is established
for the covariant and contravariant tensor irreps.

\normalsize  
\pagebreak 
\section*{Quantum supergroup $GL_q(m|n)$}
We will work on the complex field ${\bf C}$. For conveninece, 
the Lie superalgebra $gl(m|n)$ will be denoted by ${\bf g}$, and 
its quantized universal enveloping algebra by $U_q({\bf g})$. 
Here $U_q({\bf g})$ is of  
Jimbo type with $q$ being specialized to a generic complex 
parameter.   We will largely follow the notations of 
\cite{I} and \cite{II}. 
Now the  generators of $U_q({\bf g})$ are  
$K_a,  \ K_a^{-1}, \  E_{b\   {b+1}},$ $ \ E_{b+1,   b}, 
\ a\in {\bf I}, 
\ b\in {\bf I}'$,
where ${\bf I}=\{1, 2, ..., m+n\}$, 
${\bf I}'={\bf I}\backslash \{m+n\}$, 
and the defining relations are given in ~\cite{I}. 
Corresponding to any subset $\bf\Theta$ of ${\bf I^\prime}$, 
we introduce $U_q({\bf k})$ which is 
generated by the elements of  
${\cal S}_k= \{ K_a^{\pm 1},\  E_{c\, c+1}, \ E_{c+1\, c},\ 
a\in {\bf I},  \ c\in {\bf\Theta} \}$, and 
$U_q({\bf p}_\pm)$ respectively generated by those  
of ${\cal S}_k \cup \{ E_{c\ c+1},
 \ c\in {\bf I^\prime} \backslash {\bf\Theta}\}$
and ${\cal S}_k \cup \{ E_{c+1\ c}, 
\ c\in {\bf I^\prime} \backslash {\bf\Theta}\}$.
Clearly, $U_q({\bf p}_\pm)$ are ${\bf Z}_2$ - graded 
Hopf subalgebras of $U_q({\bf g})$, which will  
be called parabolic. Also, $U_q({\bf k})$ is a 
${\bf Z}_2$ - graded Hopf subalgebra of $U_q({\bf p}_\pm)$.

The representation theory of $U_q({\bf g})$ 
has been systematically developed~\cite{I}. 
Let $\pi$ be the contravariant vector representation of $U_q({\bf g})$ 
afforded by the irreducible module $E$ in the standard basis 
$\{ v_a\}$ such that $E_{a\, a\pm 1} v_b = \delta_{b\, a\pm 1} v_a$, 
$K_a v_b = q_a^{\delta_{a b}} v_b$. 
Although finite dimensional representations of $U_q({\bf g})$ are 
not completely reducible in general, 
repeated tensor products of $\pi$ can be reduced into 
direct sums of contravariant tensor irreps. Let $\Lambda^{(1)}$ 
be the set of the highest weights of all such irreps.
Let ${\overline E}$ be the dual module of $E$, and  
${\overline \pi}$ the covariant vector representation 
(dual to $\pi$).  
Repeated tensor products of ${\overline \pi}$ are again 
completely reducible. 
We denote by $\Lambda^{(2)}$ the highest weights of all the covariant
tensor irreps. (The trivial representation is 
tensorial by convention.)  
$\Lambda^{(1)}$ and $\Lambda^{(2)}$ can be explicitly 
characterized~\cite{III}. 
In particular, it is known that $\Lambda^{(1)}
\cap \Lambda^{(2)}=\{0\}$.  Also, if $\lambda\in \Lambda^{(1)}$, 
we denote by ${\bar\lambda}$ the lowest weight of the irrep with 
highest weight $\lambda$, and set $\lambda^\dagger = - {\bar\lambda}$. 
Then $\Lambda^{(2)}=\{\lambda^\dagger| 
\lambda\in \Lambda^{(1)}\}$.

It follows from standard Hopf algebra theory that $(U_q({\bf g}))^0$ 
is a Hopf superalgebra with a structure dualizing that 
of $U_q({\bf g})$.  Consider 
$T_q\subset(U_q({\bf g}))^0$ generated by 
the matrix elements of $\pi$: 
$t_{a b}\in (U_q({\bf g}))^0$, $a, b\in{\bf I}$.  
It is a ${\bf Z}_2$ - graded bi - subalgebra  
of $(U_q({\bf g}))^0$ with the generators satisfying 
an `RTT' relation~\cite{II}, but does not admit an antipode. 
Our earlier  discussions on  representations of $U_q({\bf g})$ 
imply that the matrix elements of the contravariant tensor irreps 
form a  Peter - Weyl basis for $T_q$. 
Similarly, the matrix elements ${\overline t}_{a b}\in(U_q({\bf g}))^0$ of
${\overline\pi}$ also generate a ${\bf Z}_2$ - graded bi - subalgebra
${\overline T}_q$ of $(U_q({\bf g}))^0$, for which 
the matrix elements of the covariant tensor irreps 
form a Peter - Weyl basis.

We define the algebra of functions $G_q$ on $GL_q(m|n)$ to be the 
subalgebra of $(U_q({\bf g}))^0$ generated by $\{t_{a b}, \ 
{\overline t}_{a b} | a, b\in{\bf I}\}$. It inherits a bi - superalgebra 
structure from $T_q$ and ${\overline T}_q$, and also admits an antipode. 
Thus $G_q$  is a Hopf superalgebra. 
By considering the universal $R$ - matrix of $U_q({\bf g})$ we can 
easily show that the following relation is satisfied in $G_q$
\bea 
R^{{\overline\pi}\, \pi}_{1 2}\ {\overline t}_1 \ t_2
&=& t_2 \ {\overline t}_1 \ R^{{\overline\pi}\, \pi}_{1 2}, \label{mix} 
\eea 
where $R^{{\overline\pi}\, \pi} =
q^{-\sum_{a\in{\bf I}} e_{a\, a}\otimes e_{a\, a} (-1)^{[a]} }
- (q-q^{-1})\sum_{a<b} e_{b\, a}\otimes e_{b\, a}
(-1)^{[a]+ [b] + [a][b]}$.
An immediate consequence of (\ref{mix}) is the factorization $G_q
=T_q {\overline T}_q$.
It can also be shown~\cite{II} that $G_q$ 
separates points of $U_q({\bf g})$,
that is, for any nonvanishing $x\in U_q({\bf g})$,  
there exists $f\in G_q$ such that $f( x )\ne 0$. 

\section*{Induced representations and quantum supervector bundles}
Let us introduce two left actions $L$ and $R$ 
of $U_q({\bf g})$ on $G_q$. For any $x\in U_q({\bf g})$ 
and $f\in G_q$, we define 
$R_x(f)=\sum_{(f)} (-1)^{[x]([f]+[x])}\ f_{(1)} \langle f_{(2)},
\ x\rangle$ and   
$L_x(f) =\sum_{(f)} \langle f_{(1)}, \ S^{-1}(x) \rangle f_{(2)}$,
where Sweedler's sigma notation is employed.   
We {\em trivially} extend them to $V\otimes_{\bf C} G_q$, 
and still denote the resultant actions by the same notations. 
Clearly, $L$ and $R$ graded - commutate. Now define  
\bea 
{\cal A}_q&:=&\{ f\in G_q | R_y(f)=\epsilon(y) f, \ \forall 
y\in U_q({\bf k})\},   
\eea 
which forms a subalgebra of $G_q$. 
It can be  regarded as  the superalgebra of functions on a 
quantum homogeneous superspace. 
Given any finite dimensional $U_q({\bf k})$ - module $V$, we define
\bea  
{\cal E}(V) &:=& \{ \zeta\in V\otimes_{\bf C} G_q \, |\,
R_y(\zeta)  = (\, S(y) \otimes id_{G_q}) \zeta,  \ \forall
\ y\in U_q({\bf k}) \}.  
\eea 
Then it immediately follows from  the commutativity of $R$ and $L$ 
that 
${\cal E}(V)$ furnishes a left $U_q({\bf g})$ - module under $L$. 
We will call ${\cal E}(V)$ an induced $U_q({\bf g})$ - module. 
There exists the following quantum Frobenius reciprocity: 
{\em if $W$ is a quotient $U_q({\bf g})$ - module of 
$\bigoplus_{k, l}\left\{ E^{\otimes k}\otimes
{\overline E}^{\otimes l}\right\}$, 
then there is the canonical isomorphism
\bea
\mbox{Hom}_{U_q({\bf g})} ( W, \ {\cal E}(V) )
&\cong& \mbox{Hom}_{U_q({\bf k})} ( W, \ V). 
\eea} 

To understand the underlying geometry of induced representations, 
observe that ${\cal E}(V)$ forms a two - sided ${\cal A}_q$ 
- module under the multiplication of $G_q$.  If $W$ is a 
$U_q({\bf g})$ - module, then ${\cal E}(W)$ is free.  
%i.e., ${\cal E}(W)\cong W\otimes_{\bf C} {\cal A}_q$, both as 
%a left and a right ${\cal A}_q$ - module, although the 
%isomorphisms in the two cases are different. 
It follows that  if $V$ is a  
$U_q({\bf k})$ - module which is contained as a direct summand  
in a $U_q({\bf g})$ - module, then ${\cal E}(V)$ is project, 
and in this case, we may regard ${\cal E}(V)$ as the space 
of sections of a quantum supervector bundle~\cite{Connes}
associated with ${\cal A}_q$. 
When $m\not\in{\bf\Theta}$, 
every finite dimensional $U_q({\bf k})$ - module 
$V$ has this property. Also, for general ${\bf\Theta}$,  
${\cal E}(V)$ is  projective if it 
yields non - zero ${\cal O}(\mu)$ or  
${\overline{\cal O}}(\mu)$(defined below).

Let $V$ be a finite dimensional
irreducible module over $U_q({\bf p})=U_q({\bf p}_\pm)$ 
with the  $U_q({\bf k})$ highest  weight $\mu$ and lowest  
weight $\tilde\mu$.  Define the subspaces 
${\cal O}(\mu)$ and ${\overline{\cal O}}(\mu)$ 
of ${\cal E}(V)$ by 
\bea
{\overline{\cal O}} (\mu) := 
     \{\zeta\in {\cal E}(V) \cap \left( V\otimes T_q\right) 
      |    R_y(\zeta)  = (S(y) \otimes id_{G_q}) \zeta,  \ \forall
\ y\in U_q({\bf p}) \};  \nonumber\\
{\cal O} (\mu) := 
\{ \zeta\in {\cal E}(V) \cap \left(V\otimes {\overline T}_q\right) 
|    R_y(\zeta)  = (S(y) \otimes id_{G_q}) \zeta,  \ \forall
\ y\in U_q({\bf p}) \}. 
\eea
Then we have the following quantum Borel - Weil theorem for  
the covariant and contravariant tensor irreps~\cite{II}:  
{\em As  $U_q({\bf g})$ - modules,
\bea
{\cal O} (\mu)&\cong&
 \left\{ \begin{array}{l l l }
  W((-{\tilde\mu})^\dagger), & \mbox{if} \ {\tilde\mu}\in -\Lambda^{(2)},
           & U_q({\bf p}) = U_q({\bf p}_+), \\
  W(\mu), & \mbox{if} \ \mu\in \Lambda^{(1)},
           & U_q({\bf p}) = U_q({\bf p}_-),\\ 
   \{0\},  &\mbox{otherwise}.
                     \end{array}\right. \nonumber \\
{\overline{\cal O}}(\mu)&\cong&
 \left\{ \begin{array}{l l l }
  W((-{\tilde\mu})^\dagger), & \mbox{if} \ {\tilde\mu}\in -\Lambda^{(1)},
           & U_q({\bf p}) = U_q({\bf p}_+), \\
 W(\mu), & \mbox{if} \ \mu\in \Lambda^{(2)},
           & U_q({\bf p}) = U_q({\bf p}_-),\\
   \{0\},  &\mbox{otherwise}, 
                     \end{array}\right.  
\eea
where the notation $W(\lambda)$  signifies the irreducible
$U_q({\bf g})$ - module with highest weight $\lambda$}.

\end{document}